# Monotone Difference Schemes for Convection-Dominated Diffusion-Reaction Equations Based on Quadratic Spline

Stelia O., Potapenko L., Sirenko I.

Faculty of Computer Sciences and Cybernetics,
Taras Shevchenko National University of Kyiv, Kyiv, Ukraine
`oleg.stelya@gmail.com, lpotapenko@ukr.net,`
`i.sirenko@gmail.com`

**Abstract.** A three-point monotone difference scheme is proposed for solving a one-dimensional non-stationary convection-diffusion-reaction equation with variable coefficients. The scheme is based on a parabolic spline and allows to linearly reproduce the numerical solution of the boundary value problem over the integral segment in the form of the function which continuous with its first derivative. The constructed difference scheme give a highly effective tool for solving problems with a small parameter at the older derivative in a wide range of output data of the problem. In the test case, numerical and exact solutions of the problem are compared with the significant dominance of the convective term of the equation over the diffusion. Numerous calculations showed the high efficiency of the new monotonous scheme developed.

## 1 Introduction

A large number of applied applications, in particular in ballistics, environmental modeling, medicine, etc. is based on boundary value problems for the equation of convection-diffusion-reaction type. Physical aspects of processes describing equations and analytic solutions for some partial cases can be found in [1]. Such a wide application of the equation necessitates the creation of new highly effective numerical methods for its solution, especially with the significant dominance of convection over diffusion.

Considerable attention in scientific publications is devoted to the construction of monotone difference schemes. The property of monotone schemes is that, with their numerical implementation, the monotony of the solution of the original problem with the same direction of growth is preserved. For monotonous schemes it is easy to justify their stability. In [2], ideas and approaches to numerical simulation of the interaction of convective and diffusion processes are presented. An overview of the network methods used to solve the convection-diffusion equation and examples of their application is presented.

The main problems arising in the search for approximate solutions of the equation of convection-diffusion by numerical methods are described in [3]. Finite-difference method and finite element method are considered. Various requirements for numerical algorithms are considered in [4], such as their adaptation to the peculiarities of simulated phenomena. The use of irregular and moving nets, monotonous circuits is considered. Numerical methods for solving the convection-diffusion-reaction equation are devoted to work [5, 6]. The book [7] discusses mathematical models that consist of convection-diffusion-reaction equations, written in integral, differential or weak form. In particular, qualitative properties of exact solutions of model problems of elliptic, hyperbolic and parabolic types are discussed. Some approaches to solving the convection-diffusion equation with prevailing convection, in particular the construction of difference grids adapting to the features of the solution of problems, are given in publications [8, 11, 12]. In [13,14] non-standard finite difference schemes (NSFD) are proposed for approximating solutions of the class of generalized convection-diffusion-reaction equations.

The three-point monotone difference scheme for solving a one-dimensional nonstationary convection-diffusion-reaction equation with variable coefficients, based on a parabolic spline, proposed in this work.

## 2    Formulation of the problem

We consider the boundary value problem for the nonstationary equation of convection-diffusion-reaction

$$\frac{\partial u(x,t)}{\partial t} = D(x,t)\frac{\partial^2 u(x,t)}{\partial x^2} - V(x,t)\frac{\partial u(x,t)}{\partial x} +$$

$$+ A(x,t)u(x,t) + f(x,t), \ x \in (0,L), t > 0, \qquad (1)$$

$$u(0,t) = U_0(t), \qquad (2)$$

$$u(L,t) = U_L(t), \qquad (3)$$

$$u(x,0) = g(x), \qquad (4)$$

where $V(x,t) \neq 0$, $0 < D(x,t) << 1$.

We introduce a uniform grid on the time variable

$$\Delta_t : t_{k+1} = t_k + \rho, \ k = 0,1,2,...,t_0 = 0, \rho = const.$$

Then we write the equation with the sampled time variable:

$$D(x,t_{k+1})\frac{d^2u(x,t_{k+1})}{dx^2} - V(x,t_{k+1})\frac{du(x,t_{k+1})}{dx} - \frac{1}{\rho}u(x,t_{k+1}) + A(x,t_{k+1})u(x,t_{k+1}) = \quad (5)$$

$$= -f(x,t_{k+1}) - \frac{1}{\rho}u(x,t_k),$$

$$u(0,t_{k+1}) = U_0(t_{k+1}), \quad (6)$$

$$u(L,t_{k+1}) = U_L(t_{k+1}), \quad (7)$$

$$u(x,t_0) = g(x). \quad (8)$$

Let two splitting ($\Delta_x$ and $\Delta_\tau$) are defined on the interval $[0,L]$.

a) $\Delta_x : 0 = x_0 < x_1 < ... < x_N = L,$

b) $\Delta_\tau : 0 = \tau_0 < \tau_1 < ... < \tau_{N-1} = L,$ \quad (9)

where $x_i < \tau_i < x_{i+1}, i = \overline{1, N-2}$.

We denote $C_i$ and $\varphi_i$ the values of some network functions, respectively, on the grids (9a) and (9b), where in $\varphi_0 = C_0, \varphi_{N-1} = C_N$. We shall seek the solution of the problem in the form of a parabolic spline [9]. To do this, we write the piecewise-quadratic function $C(x)$ at the time moment $t_{k+1}, x \in [\tau_i, \tau_{i+1}]$ and we find the first and second derivatives of this function on each of these segments, substituting them together with the function itself in equation (5). We get:

$$D^{k+1}(x)\left\{\varphi_i \frac{2}{(x_{i+1}-\tau_i)(\tau_{i+1}-\tau_i)} - C_{i+1}\frac{2}{(x_{i+1}-\tau_i)(\tau_{i+1}-x_{i+1})} + \varphi_{i+1}\frac{2}{(\tau_{i+1}-x_{i+1})(\tau_{i+1}-\tau_i)}\right\} -$$

$$-V^{k+1}(x)\left\{\varphi_i \frac{(x-x_{i+1})+(x-\tau_{i+1})}{(x_{i+1}-\tau_i)(\tau_{i+1}-\tau_i)} - C_{i+1}\frac{(x-\tau_i)+(x-\tau_{i+1})}{(x_{i+1}-\tau_i)(\tau_{i+1}-x_{i+1})} + \varphi_{i+1}\frac{(x-\tau_i)+(x-x_{i+1})}{(\tau_{i+1}-x_{i+1})(\tau_{i+1}-\tau_i)}\right\} +$$

$$+\left(A^{k+1}(x)-\frac{1}{\rho}\right)\left\{\varphi_i \frac{(x-x_{i+1})(x-\tau_{i+1})}{(x_{i+1}-\tau_i)(\tau_{i+1}-\tau_i)} - C_{i+1}\frac{(x-\tau_i)(x-\tau_{i+1})}{(x_{i+1}-\tau_i)(\tau_{i+1}-x_{i+1})} +\right.$$

$$\left.+\varphi_{i+1}\frac{(x-\tau_i)(x-x_{i+1})}{(\tau_{i+1}-x_{i+1})(\tau_{i+1}-\tau_i)}\right\} = -f^{k+1}(x) - \frac{1}{\rho}u^k(x). \quad (10)$$

Without limiting the generality of the algorithm, let's put $D = const, V = const, A = const,$

Also, let's put $x = x_{i+1}$. We get:

$$\varphi_i \frac{2D}{(x_{i+1}-\tau_i)(\tau_{i+1}-\tau_i)} - C_{i+1}\frac{2D}{(x_{i+1}-\tau_i)(\tau_{i+1}-x_{i+1})} + \varphi_{i+1}\frac{2D}{(\tau_{i+1}-x_{i+1})(\tau_{i+1}-\tau_i)} -$$

$$-\varphi_i \frac{V(x_{i+1}-\tau_{i+1})}{(x_{i+1}-\tau_i)(\tau_{i+1}-\tau_i)} + C_{i+1}\frac{V((x_{i+1}-\tau_i)+(\tau_{i+1}-x_{i+1}))}{(x_{i+1}-\tau_i)(\tau_{i+1}-x_{i+1})} - \varphi_{i+1}\frac{V(x_{i+1}-\tau_i)}{(\tau_{i+1}-x_{i+1})(\tau_{i+1}-\tau_i)} +$$

$$+\left(A-\frac{1}{\rho}\right)C_{i+1} = -f^{k+1}(x_{i+1}) - \frac{1}{\rho}u^k(x_{i+1}).$$

From the last expression we will define $C_{i+1}$

$$C_{i+1} = \left\{\varphi_i\left[\frac{2D-V(x_{i+1}-\tau_{i+1})}{(x_{i+1}-\tau_i)(\tau_{i+1}-\tau_i)}\right] + \varphi_{i+1}\left[\frac{2D-V(x_{i+1}-\tau_i)}{(\tau_{i+1}-x_{i+1})(\tau_{i-1}-\tau_i)}\right] + f_{i+1} + \left(\frac{1}{\rho}+A\right)u^k_{i+1}\right\}:$$

$$:\left[\frac{2D-V((x_{i+1}-\tau_i)+(x_{i+1}-\tau_{i+1}))}{(x_{i+1}-\tau_i)(\tau_{i+1}-x_{i+1})} + \frac{1}{\rho} - A\right].$$

Similarly, for a segment $[\tau_{i-1}, \tau_i]$ we put $x = x_i$ and we find $C_i$. Then, from the conditions of continuity of the first derivatives of functions $C(x)$ in points $\tau_i$ and taking into account the found values $C_i$ and $C_{i+1}$, we obtain

$$\varphi_{i-1}Q_{i-1} - \varphi_i Q_i + \varphi_{i+1}Q_{i+1} = F_i, i = \overline{1, N-2},$$

where

$$Q_{i-1} = \frac{\dfrac{(\tau_i-x_i)}{(x_i-\tau_{i-1})(\tau_i-\tau_{i-1})} - \left[\dfrac{(2D-V(x_i-\tau_i))(\tau_i-\tau_{i-1})}{(x_i-\tau_{i-1})(\tau_i-\tau_{i-1})(x_i-\tau_{i-1})(\tau_i-x_i)}\right]}{\left[\dfrac{2D-V((x_i-\tau_{i-1})+(x_i-\tau_i))}{(x_i-\tau_{i-1})(\tau_i-x_i)} + \dfrac{1}{\rho} - A\right]},$$

$$Q_i = \left\{\frac{\left[\dfrac{(2D-V(x_i-\tau_{i-1}))(\tau_i-\tau_{i-1})}{(\tau_i-x_i)(\tau_i-\tau_{i-1})(x_i-\tau_{i-1})(\tau_i-x_i)}\right]}{\left[\dfrac{2D-V((x_i-\tau_{i-1})+(x_i-\tau_i))}{(x_i-\tau_{i-1})(\tau_i-x_i)} + \dfrac{1}{\rho}-A\right]} - \right.$$

$$\left. -\frac{(\tau_i-\tau_{i-1})+(\tau_i-x_i)}{(\tau_i-x_i)(\tau_i-\tau_{i-1})} + \frac{(\tau_i-x_{i+1})+(\tau_i-\tau_{i+1})}{(x_{i+1}-\tau_i)(\tau_{i+1}-\tau_i)} - \right.$$

$$-\left[\frac{(2D-V(x_{i+1}-\tau_{i+1}))}{(x_{i+1}-\tau_i)(\tau_{i+1}-\tau_i)}\right] \times \frac{\left[\frac{(\tau_i-\tau_{i+1})}{(x_{i+1}-\tau_i)(\tau_{i+1}-x_{i+1})}\right]}{\left[\frac{2D-V((x_{i+1}-\tau_i)+(x_{i+1}-\tau_{i+1}))}{(x_{i+1}-\tau_i)(\tau_{i+1}-x_{i+1})}+\frac{1}{\rho}-A\right]}\right\},$$

$$Q_{i+1} = \frac{\frac{[2D-V(x_{i+1}-\tau_i)](\tau_i-\tau_{i+1})}{(\tau_{i+1}-x_{i+1})(\tau_{i+1}-\tau_i)(x_{i+1}-\tau_i)(\tau_{i+1}-x_{i+1})}}{\left[\frac{2D-V((x_{i+1}-\tau_i)+(x_{i+1}-\tau_{i+1}))}{(x_{i+1}-\tau_i)(\tau_{i+1}-x_{i+1})}+\frac{1}{\rho}-A\right]} - \frac{(\tau_i-x_{i+1})}{(\tau_{i+1}-x_{i+1})(\tau_{i+1}-\tau_i)}\right\},$$

$$F_i = -\frac{\left[f_{i+1}+\left(\frac{1}{\rho}-A\right)u^k_{i+1}\right](\tau_i-\tau_{i+1})}{(x_{i+1}-\tau_i)(\tau_{i+1}-x_{i+1})} + \left[\frac{2D-V((x_{i+1}-\tau_i)+(x_{i+1}-\tau_{i+1}))}{(x_{i+1}-\tau_i)(\tau_{i+1}-x_{i+1})}+\frac{1}{\rho}-A\right]$$

$$+\left\{\frac{\left[f_i+\left(\frac{1}{\rho}-A\right)u^k_i\right](\tau_i-\tau_{i-1})}{(x_i-\tau_{i-1})(\tau_i-x_i)}\right\} : \left[\frac{2D-V((x_i-\tau_{i-1})+(x_i-\tau_i))}{(x_i-\tau_{i-1})(\tau_i-x_i)}+\frac{1}{\rho}-A\right].$$

Note that from the time $k=1$ as values $u(x_{i+1})$ and $u(x_i)$ we use the values of $C_{i+1}$ and $C_i$. We consider a uniform grid $\Delta_\tau$ with step $h$. The grid $\Delta_x$ is placed relative to the grid $\Delta_\tau$ so that $\tau_i - x_i = \mu$. After the expression transformations of the expressions $Q_{i-1}$, $Q_i$, $Q_{i+1}$ and $F_i$, we obtain a system of equations for the internal nodes of the grid at the time $t_{k+1}$

$$\alpha\varphi^{k+1}_{\tau_{i-1}} - \gamma\varphi^{k+1}_{\tau_i} + \beta\varphi^{k+1}_{\tau_{i+1}} = -\frac{f^{k+1}_{x_i} + f^{k+1}_{x_{i+1}}}{2} - \frac{u^k_{x_{i+1}} + u^k_{x_i}}{2\rho}, \quad i = \overline{1, N-2}, \tag{11}$$

where

$$\alpha = a + \frac{\mu^2}{2h^2}\left(\frac{1}{\rho}-A\right), \quad \gamma = a+b+\frac{2h^2-\mu^2-(h-\mu)^2}{2h^2}\left(\frac{1}{\rho}-A\right),$$

$$\beta = b - \frac{(h-\mu)^2}{2h^2}\left(\frac{1}{\rho}-A\right), \quad a = \frac{D}{h^2}+\frac{V}{h^2}\mu, \quad b = \frac{D}{h^2}-\frac{V}{h}\left(1-\frac{\mu}{h}\right).$$

**Theorem.** The difference scheme (11) is monotone under the condition $\rho > \max(\rho_1, \rho_2)$ where

$$\rho_1 = \frac{\mu^2}{2(D+V\mu)+A\mu^2}, \quad \rho_2 = \frac{(h-\mu)^2}{2(D+V\mu-Vh)+A(h-\mu)^2}$$

and the condition $\rho \leq \frac{1}{A}$, if $A > 0$.

*Proving.* From the condition $\rho > \rho_1$ we have $\alpha > 0$. From the condition $\rho > \rho_2$ we have $\beta > 0$. Then

$$\gamma - \alpha - \beta = \frac{2h^2 - \mu^2 - (h-\mu)^2}{2h^2}\left(\frac{1}{\rho} - A\right) + \frac{\mu^2}{2h^2}\left(\frac{1}{\rho} - A\right) + \frac{(h-\mu)^2}{2h^2}\left(\frac{1}{\rho} - A\right) = \frac{1}{\rho} - A \geq 0.$$

Thus $\gamma \geq \alpha + \beta$, this means that the difference scheme (11) is monotone under defined conditions [3].

The study of the error of the scheme approximation is given in [10]. Note that the presence of the reaction in equation (1) does not affect the accuracy of this scheme.

## 3    Examples of calculations

The numerical calculations were carried out with the help of the developed software.

It is easy to verify that the function $u(x,t) = \frac{1}{2\sqrt{\pi D(t+1)}} e^{-\frac{(x-0,4-t)^2}{4D(t+1)}}$ satisfies the equation (1) for $V(x,t) = 1$, $A(x,t) = 0$, $f(x,t) = 0$.

The initial and boundary conditions were given as follows:

$$u(x,0) = \frac{1}{2\sqrt{\pi D}} e^{-\frac{(x-0.4)^2}{4D}},$$

$$u(0,t) = \frac{1}{2\sqrt{\pi D(t+1)}} e^{-\frac{(-0,4-t)^2}{4D(t+1)}}, \quad u(L,t) = \frac{1}{2\sqrt{\pi D(t+1)}} e^{-\frac{(2-t)^2}{4D(t+1)}}.$$

For the numerical solution of the problem, we used such values $L = 2.4$, $D = 0.001$, $h = 0.0005$.

Comparison of the numerical solution with the exact solution is given on Fig. 1.

The numerical solution of the problem can be recreated in the form of a parabolic spline over the entire integration segment.

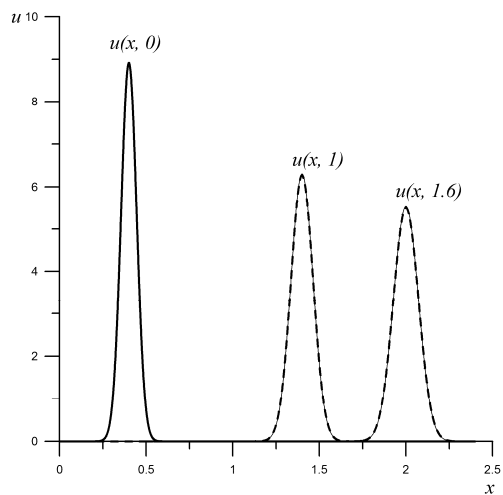

**Fig. 1.** Exact (solid line) and numerical (dotted line) solution of the problem at different times

## 4      Conclusions

The theoretical studies and the results of numerical experiments show that the proposed monotone difference scheme for the convection-diffusion-reaction equation allows solving boundary-value problems for a wide range of the coefficients values of the equation. This is especially true in cases where convection greatly exceeds diffusion. The monotone difference scheme provides the stability of the numerical solution, and the use of a parabolic spline for its construction allows you to reproduce the solution in the form of a continuous function for each time point.